\documentclass[10pt]{article}
\usepackage{hyperref}
\usepackage{amsfonts}
\usepackage{enumerate}
\usepackage{enumitem}
\usepackage{graphicx}

\usepackage{amsmath}
\usepackage{mathtools}
\usepackage{setspace}
\usepackage{color}
\usepackage{array}
\usepackage[font={small}]{caption}
\usepackage{eurosym}
\usepackage{amsmath}
\usepackage{amsfonts}
\usepackage{amssymb}
\usepackage{pifont}
\usepackage{graphics}
\usepackage{array}
\usepackage{shortvrb}
\usepackage{epsf}
\usepackage{graphicx}
\usepackage{caption}
\usepackage{rotating}
\usepackage{wasysym}
\usepackage{multirow}
\usepackage{empheq}
\usepackage{acronym}
\usepackage{epstopdf}
\usepackage{flushend}
\usepackage{caption}
\usepackage{comment}
\usepackage{float}
\usepackage{subcaption}
\usepackage{multirow,array}
\usepackage{lineno}
\usepackage{upgreek}
\usepackage{amsmath}
\usepackage{amsthm}

\graphicspath{{./figs/}}
\DeclareGraphicsExtensions{.pdf}

\begin{document}
\begin{center}
{\large\bf {Fractional Operators for Nonlinear Electrical Circuits}}

\vskip.20in

Ioannis Dassios$^{1*}$, 
\\[2mm]
{\footnotesize
$^1$Aristotle University of Thessaloniki, Thessaloniki, Greece
\\[5pt]
$^*$Corresponding author: ioannisdassios@gmail.com
\\[5pt]
}
\end{center}

{\footnotesize
\noindent
\textbf{Abstract:} 
This article introduces two new fractional operators with sine ($\sin$) and cosine ($\cos$) kernels, motivated by their fundamental role in modeling AC signals in electrical circuits. The operators are designed to improve the analysis of nonlinear components such as the memristor by transforming certain nonlinear equations into simpler linear forms, particularly in systems with memory effects. 
}

\section{Introduction}

Trigonometric functions such as sine ($\sin$) and cosine ($\cos$) play a central role in electrical engineering, particularly in the analysis of alternating current (AC) circuits and power systems~ \cite{brown_trigonometric_2022,Das4,Das6,Das7,Milano,smith_trigonometric_2023,wang_trigonometric_2021}. They are essential for describing periodic signals, as well as voltage, current, and phase behavior in circuits. Beyond classical circuit theory, trigonometric functions also appear in control systems and signal processing, including applications in robotics, telecommunications, and audio processing.

In recent years, the development of nonlinear electrical components, most notably, the memristor~ \cite{Chua1,Chua2,Chua3}, has motivated the search for advanced mathematical tools capable of capturing complex dynamics, including memory effects. The memristor, a passive two-terminal component whose resistance depends on the history of the applied voltage, is of particular interest due to its applications in neuromorphic computing and next-generation memory systems.
Fractional operators offer powerful tools for modeling complex systems, capturing memory effects and long-range interactions that are prevalent in practical applications, see \cite{Das3,Das8,Das9,Das10,Das11,Das12,Kac,fractional_operators}. Among the various tools in fractional calculus, the Caputo fractional derivative (C) is widely used in engineering due to its compatibility with classical initial conditions:
\begin{equation}\label{eq1}
\mathcal{D}_c^\alpha f(t) = \frac{1}{\Gamma(1 - \alpha)} \int_{0}^{t} (t - \tau)^{-\alpha} \frac{d}{d\tau} f(\tau) \, d\tau, \quad 0 < \alpha \leq 1,
\end{equation}
where \(\Gamma(\cdot)\) denotes the Gamma function.
Alternative formulations include the Caputo-Fabrizio (CF) operator:
\begin{equation}\label{eq2}
\mathcal{D}_{cf}^\alpha f(t) = \frac{B(\alpha)}{1 - \alpha} \int_{0}^{t} e^{-\frac{\alpha}{1 - \alpha}(t - \tau)} \frac{d}{d\tau} f(\tau) \, d\tau,
\end{equation}
and the Atangana-Baleanu (AB) operator:
\begin{equation}\label{eq3}
\mathcal{D}_{ab}^\alpha f(t) = \frac{B(\alpha)}{1 - \alpha} \int_{0}^{t} E_\alpha\left[-\alpha \frac{(t - \tau)^\alpha}{1 - \alpha}\right] \frac{d}{d\tau} f(\tau) \, d\tau,
\end{equation}
where \(B(\alpha)\) is a normalization function with \(B(0) = B(1) = 1\), and \(E_\alpha(z)\) denotes the one-parameter Mittag-Leffler function. These operators introduce exponential and non-singular kernels, providing additional flexibility in modeling nonlocal and memory-dependent phenomena~ \cite{At,Ca,Das1,Das2,Das5}. 

In this paper, we introduce two novel fractional operators with trigonometric kernels, developed to support the analysis of nonlinear electrical circuits. These operators offer an alternative mathematical framework for modeling components such as the memristor, where memory effects and nonlinear behavior are significant.

\section{Main Results}

\begin{figure}[h]
\begin{center}
\includegraphics[width=0.5\textwidth]{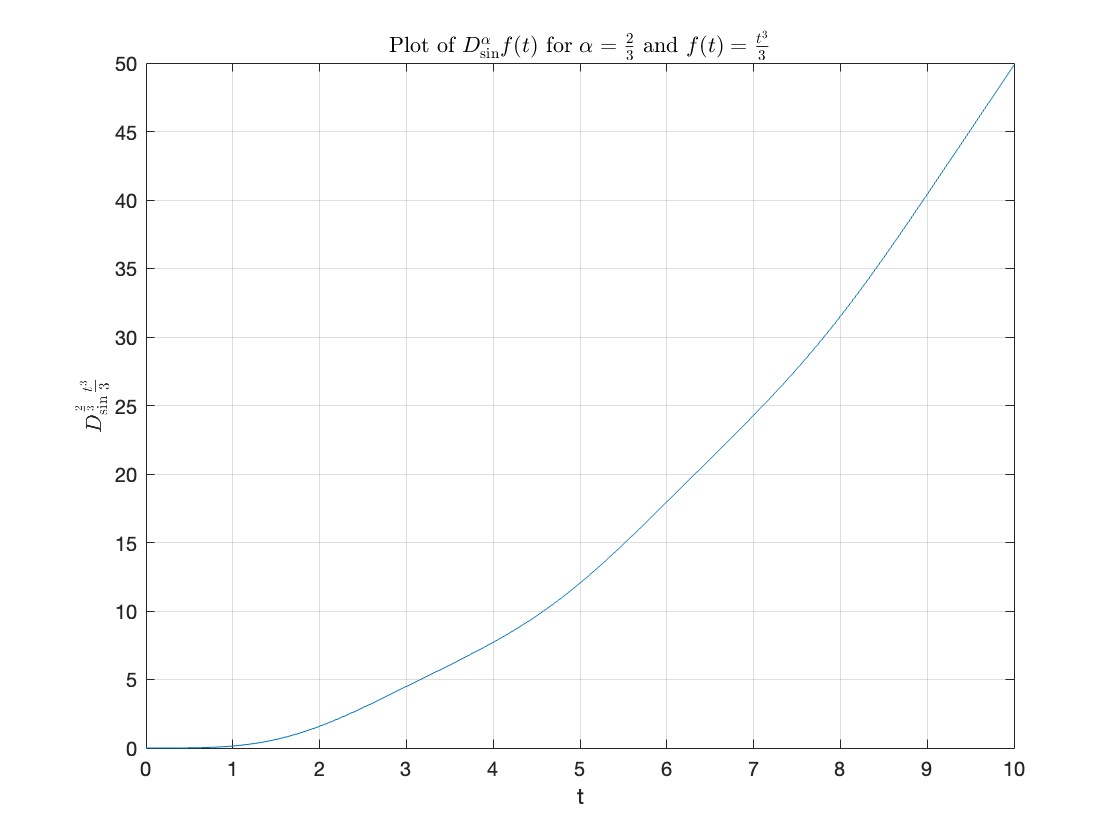}\includegraphics[width=0.5\textwidth]{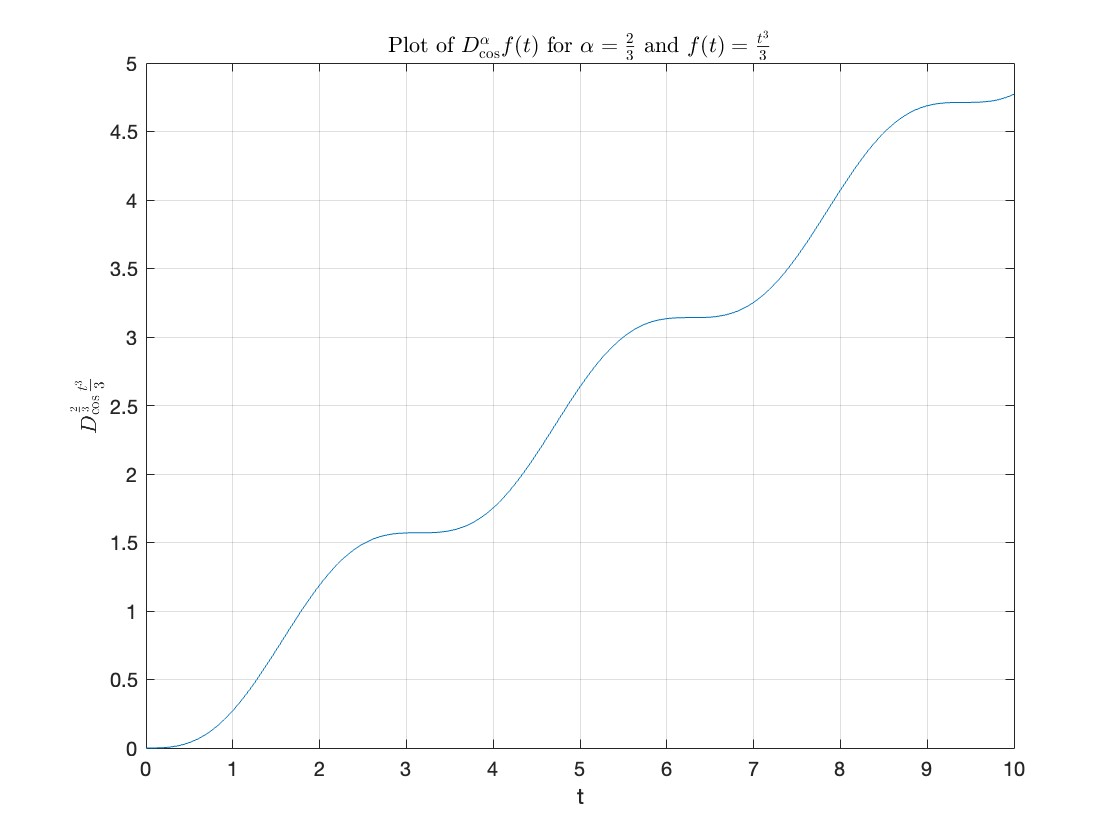}
\caption{Plots of the (DS) (left) and (DC) (right) operators for $\alpha = \frac{2}{3}$ and $f(t) = \frac{t^3}{3}$.}
\label{Figure 1}
\end{center}
\end{figure}
\begin{figure}[h]
\begin{center}
\includegraphics[width=0.50\textwidth]{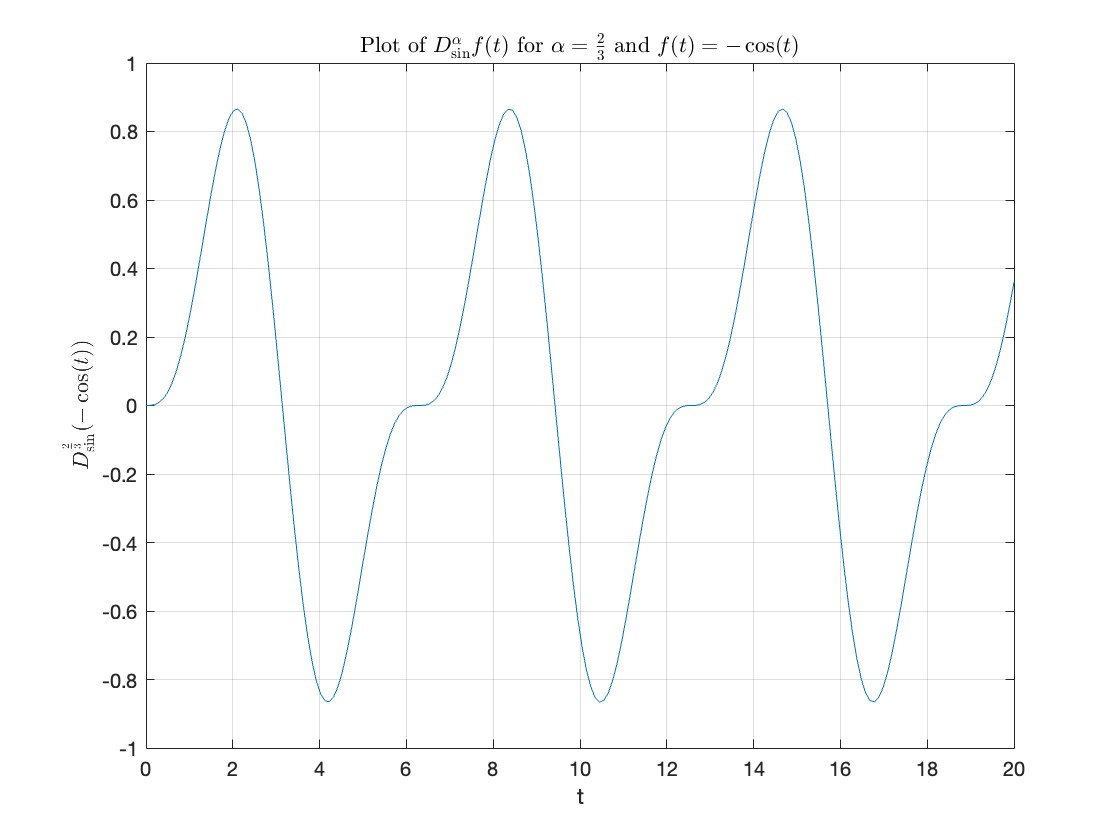}\includegraphics[width=0.50\textwidth]{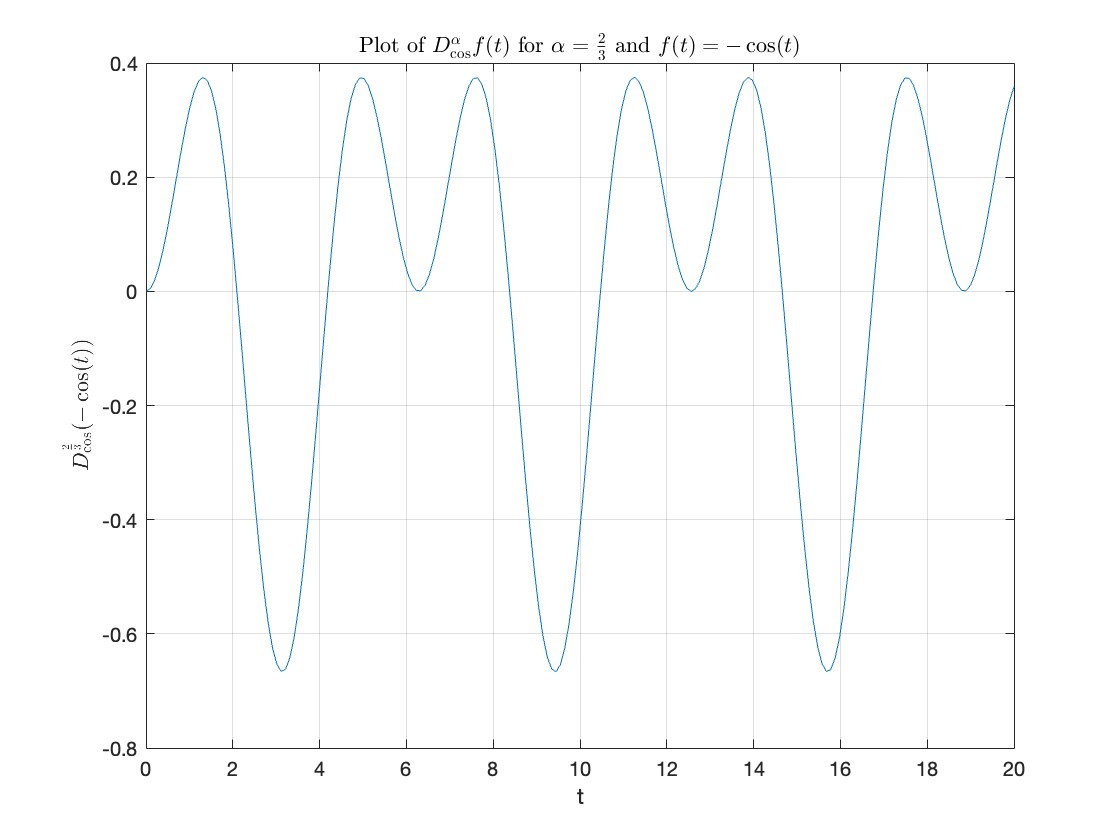}
\caption{Plots of the (DS) (left) and (DC) (right) operators for $\alpha=\frac{2}{3}$ and $f(t)=-\cos(t)$.}
\label{Figure 2}
\end{center}
\end{figure}
\begin{figure}[h]
\begin{center}
\includegraphics[width=0.50\textwidth]{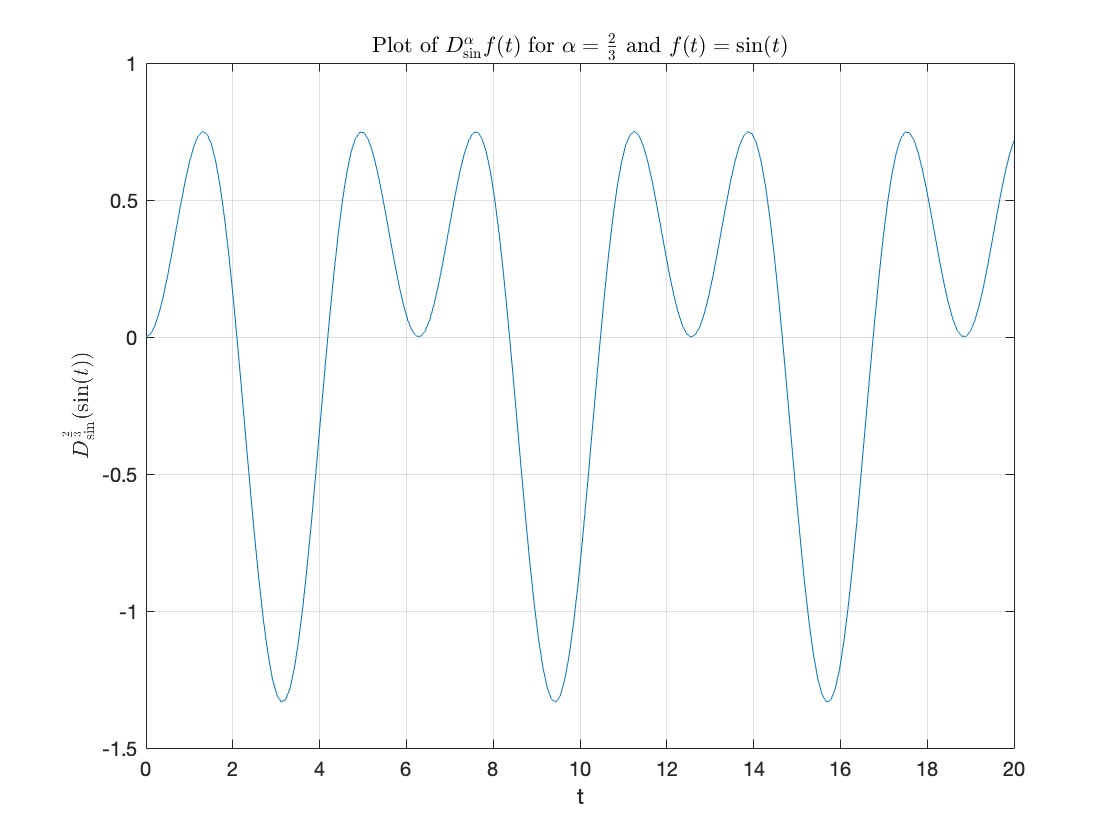}\includegraphics[width=0.50\textwidth]{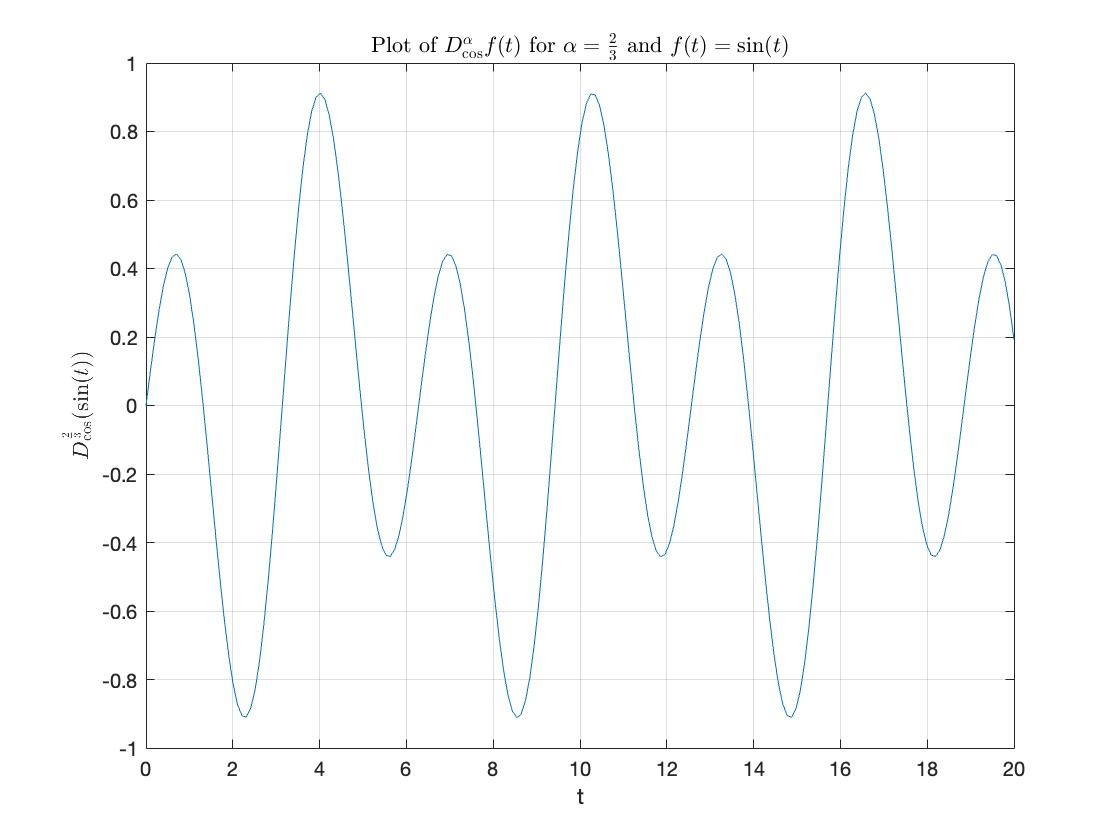}
\caption{Plots of the (DS) (left) and (DC) (right) operators for $\alpha=\frac{2}{3}$ and $f(t)=\sin(t)$.}
\label{Figure 3}
\end{center}
\end{figure}
\begin{figure}[h]
\begin{center}
\includegraphics[width=0.50\textwidth]{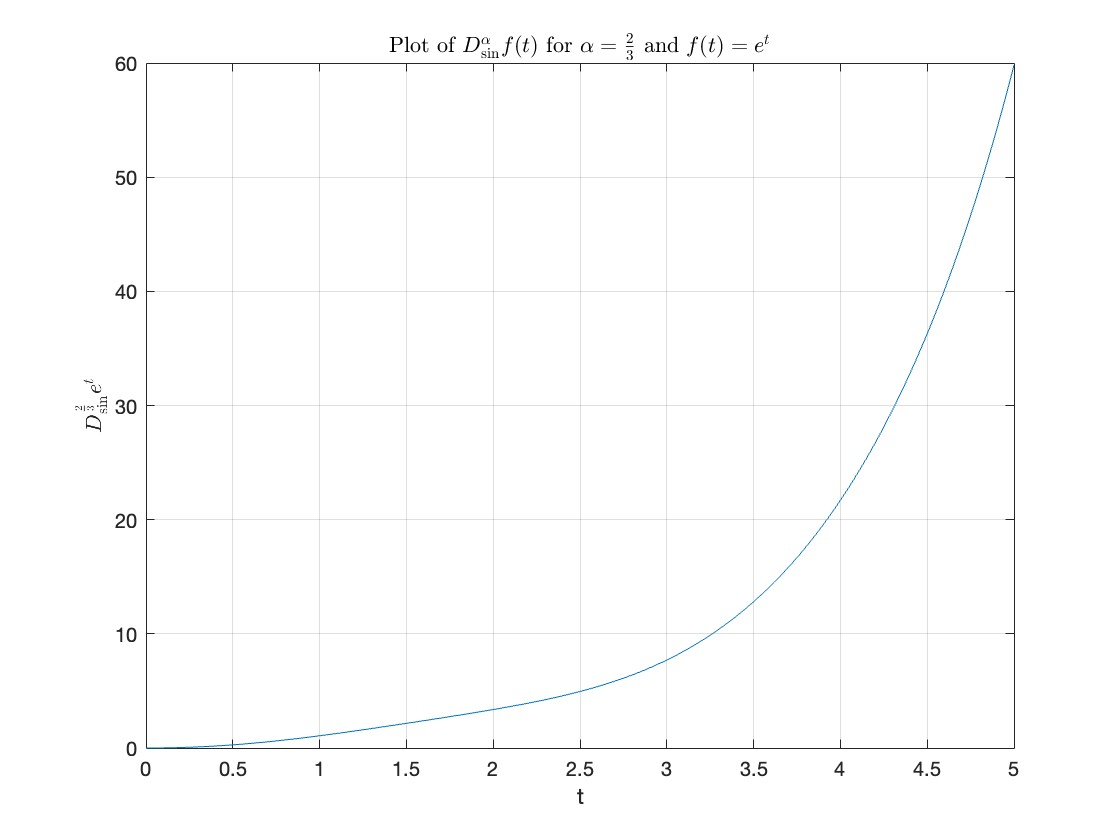}\includegraphics[width=0.50\textwidth]{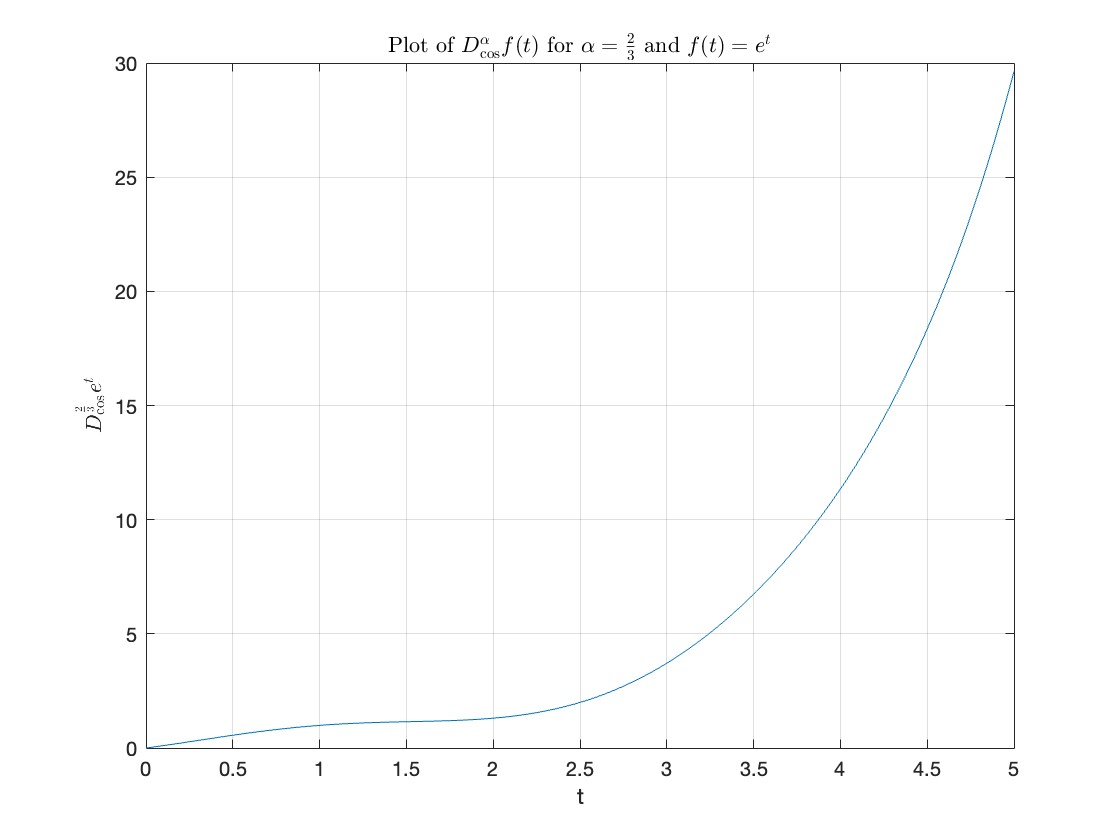}
\caption{Plots of the (DS) (left) and (DC) (right) operators for $\alpha=\frac{2}{3}$ and $f(t)=e^t$.}
\label{Figure 4}
\end{center}
\end{figure}
\begin{figure}[h]
\begin{center}
\includegraphics[width=0.49\textwidth]{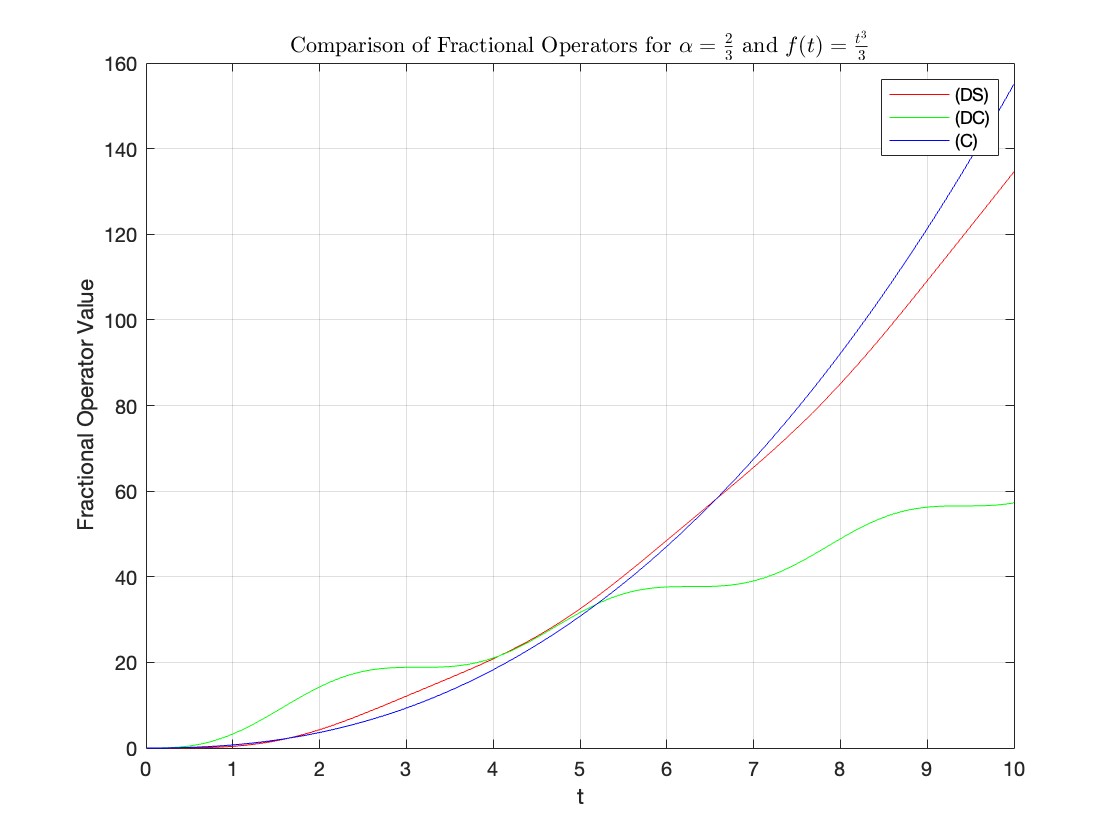}
\includegraphics[width=0.49\textwidth]{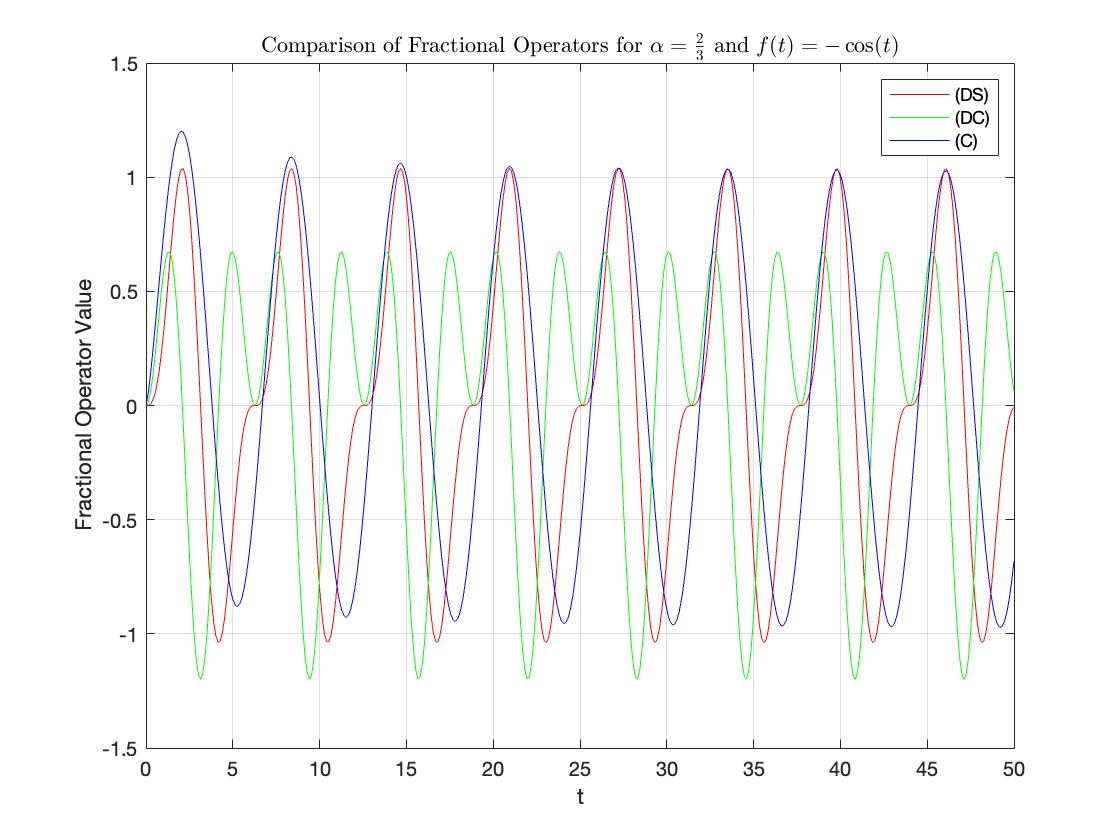}
\caption{Comparison of the (DS), (DC), and (C) operators with $\alpha = \frac{2}{3}$, applied to $f(t) = \frac{1}{3}t^3$ (left) and $f(t) = -\cos(t)$ (right).}
\label{Figure 5}
\end{center}
\end{figure}
\begin{figure}[h]
\begin{center}
\includegraphics[width=0.49\textwidth]{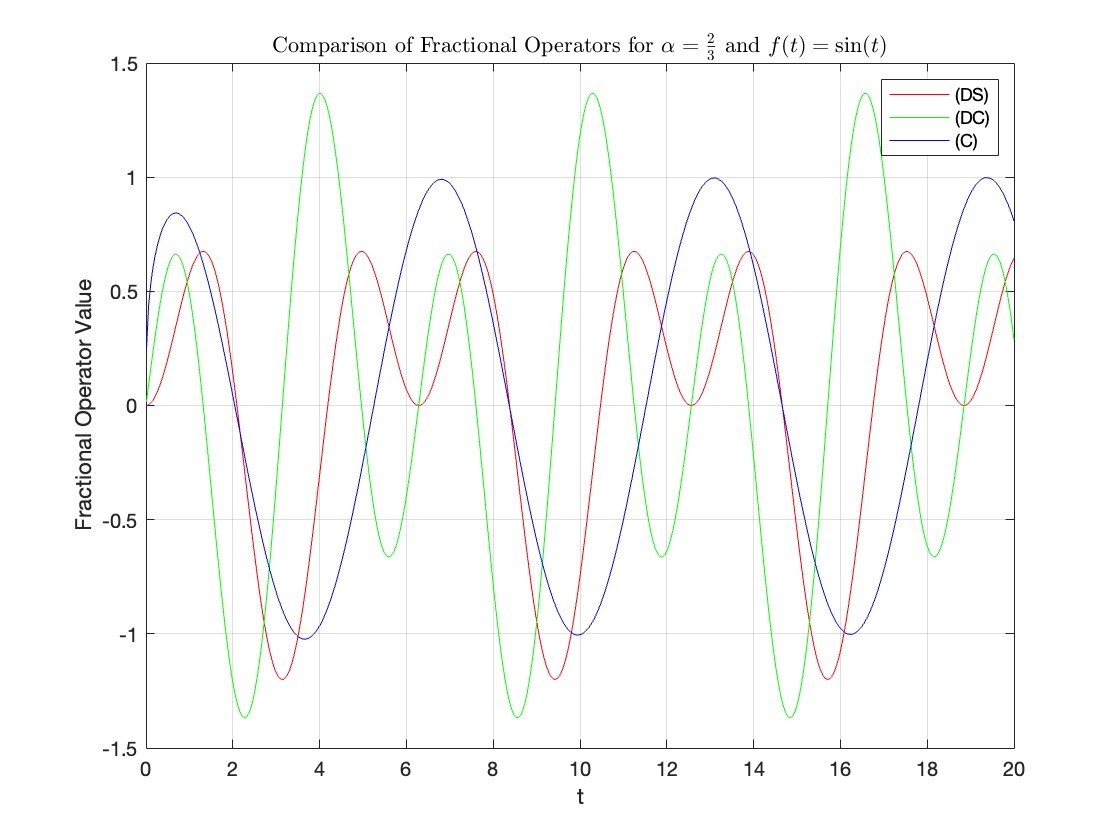}
\includegraphics[width=0.49\textwidth]{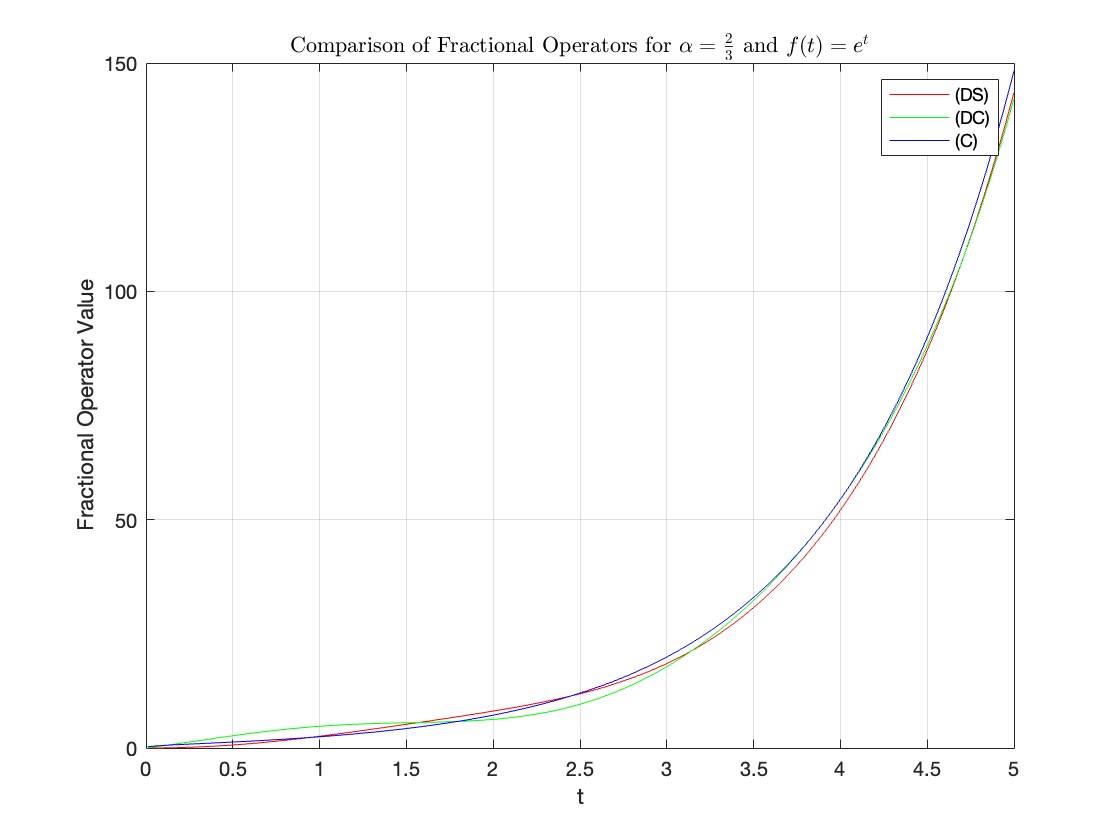}
\caption{Comparison of the (DS), (DC), and (C) operators with $\alpha = \frac{2}{3}$, applied to $f(t)=\sin(t)$ (left) and $f(t)=e^t$ (right).}
\label{Figure 6}
\end{center}
\end{figure}
\begin{figure}[h]
\begin{center}
\includegraphics[width=0.49\textwidth]{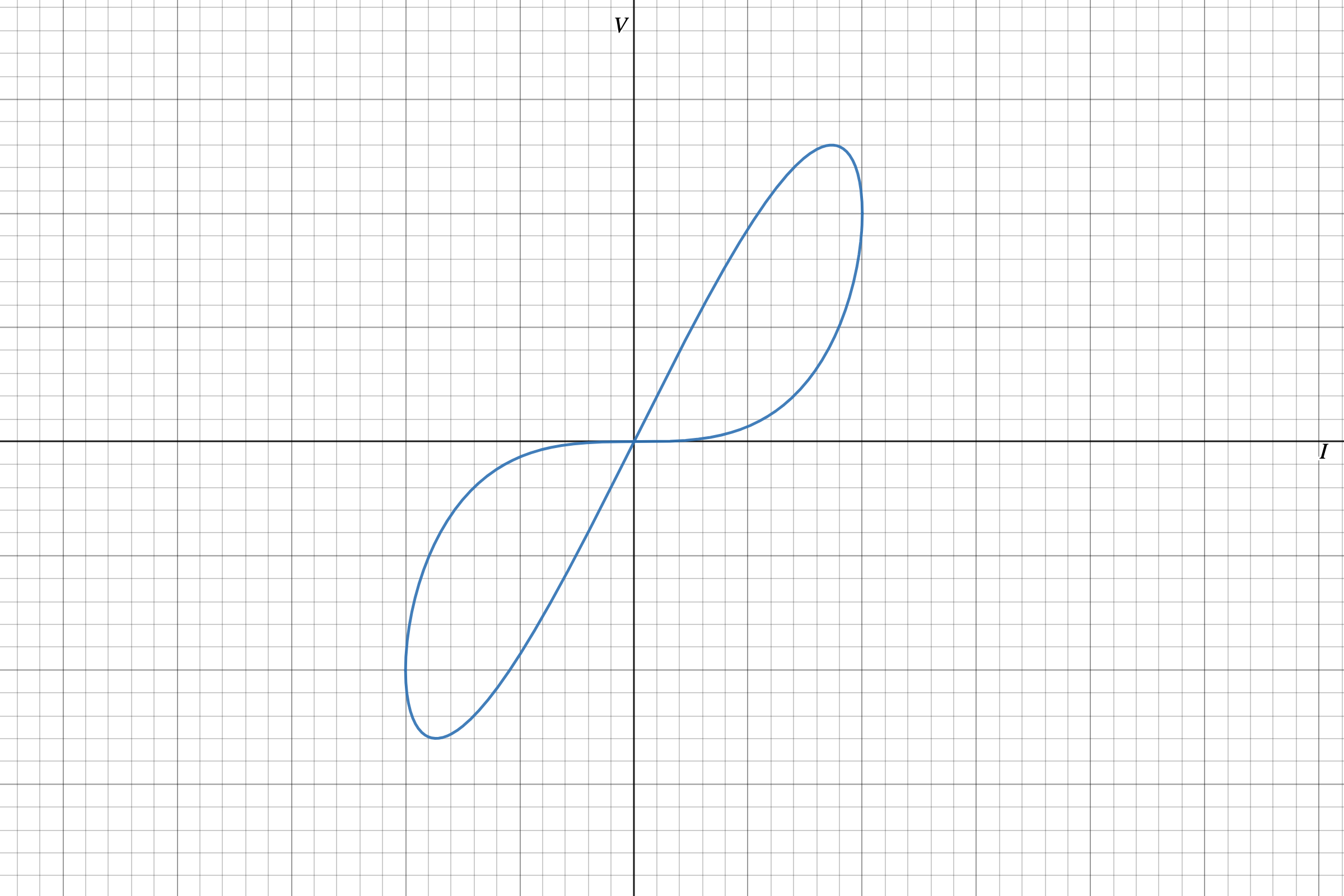}
\includegraphics[width=0.49\textwidth]{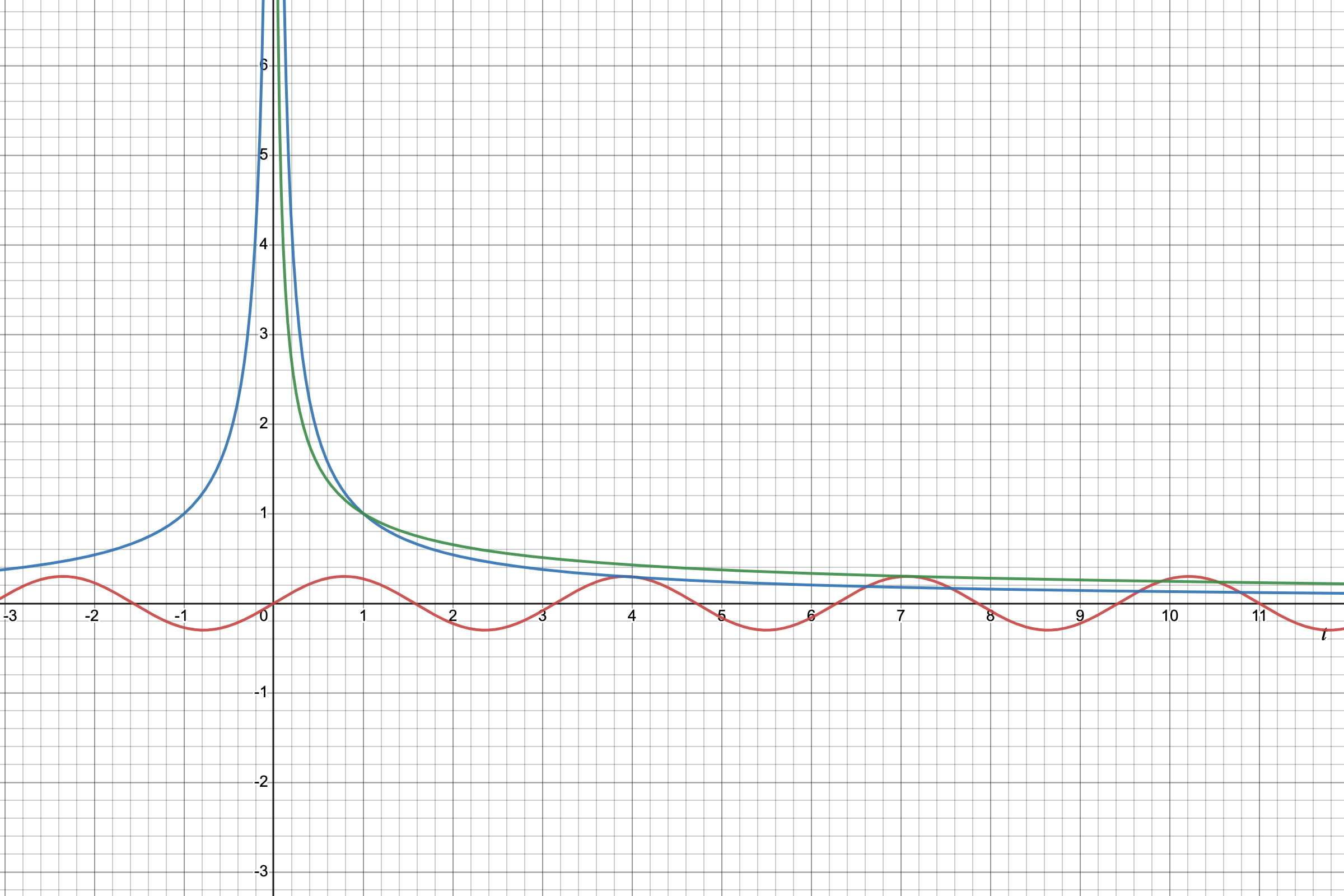}
\caption{On the left: the \(V\)-\(I\) graph of the memristor model \(V = M(q)I\). On the right: the function $\frac{3}{10} \sin(2t)$ (red), $t^{-\frac{8}{9}}$ (blue), and $t^{-\frac{43}{70}}$ (green).}
\label{Figure 7}
\end{center}
\end{figure}

In this section, we present the main contributions of the paper. The kernel associated with the (C) derivative~\eqref{eq1} is given by 
$k(t) = t^{-a}$, 
while the (CF) operator~\eqref{eq2} uses the exponential kernel 
$k(t) = e^{-\frac{a}{1 - a}t}$.
The (AB) operator~\eqref{eq3} employs the Mittag-Leffler type kernel
 $k(t)=\sum_{k=0}^\infty(-1)^k\left[\frac{a}{1-a}\right]^k\frac{t^{ak}}{\Gamma(1+ak)}$. 
We introduce two new fractional operators with trigonometric kernels defined as $k(t) = \sin\big(\frac{a}{1 - a}t\big)$ and $k(t) = \cos\big(\frac{a}{1 - a}t\big)$. These operators are developed to help model nonlinear behavior in electrical circuits more effectively. In particular, we aim for them to serve as useful alternatives in problems involving components with memory and nonlinearity, such as the memristor.
\\\\
\textbf{Definition 2.1.} We define with $\mathfrak{D}_{\sin}^\alpha$ the fractional operator (DS):
\begin{equation}\label{eq4}
\mathfrak{D}_{\sin}^\alpha f(t)
=
 \frac{N(\alpha)}{1 - \alpha} \int_{0}^{t} \sin\left[\frac{a}{1-a}(t-\tau)\right]\frac{d}{d\tau} f(\tau) \, d\tau,
\end{equation}
and with $\mathfrak{D}_{\cos}^\alpha$ the fractional operator (DC):
\begin{equation}\label{eq5}
\mathfrak{D}_{\cos}^\alpha f(t)
=
 \frac{N(\alpha)}{1 - \alpha} \int_{0}^{t} \cos\left[\frac{a}{1-a}(t-\tau)\right] \frac{d}{d\tau} f(\tau) \, d\tau,
\end{equation}
where \(N(\alpha)\) is a normalization function. 
\\\\
Let $\mathcal{L}$ be the Laplace transform with $F(s)=\mathcal{L}\left\{ f(t)\right\}$ and $f(0)$ initial condition of f(t). Then we state the following Proposition for \eqref{eq4}, \eqref{eq5}.
\\\\
\textbf{Proposition 2.1.} We consider the fractional operators \eqref{eq4}, \eqref{eq5}. Then:
\begin{equation}\label{eq6}
\mathcal{L}\left\{\mathfrak{D}_{\sin}^\alpha f(t)\right\}
=
 \frac{\alpha N(\alpha)}{(1 - \alpha)^2s^2+\alpha^2}(sF(s)-f(0)),
\end{equation}
and
\begin{equation}\label{eq7}
\mathcal{L}\left\{\mathfrak{D}_{\cos}^\alpha f(t)\right\}
=
 \frac{s(1-\alpha) N(\alpha)}{(1 - \alpha)^2s^2+\alpha^2}(sF(s)-f(0)).
\end{equation}
\textbf{Proof.} We apply the Laplace transform into \eqref{eq4} and get:
\[
\mathcal{L}\left\{\mathfrak{D}_{\sin}^\alpha f(t)\right\}
=
\mathcal{L}\left\{ \frac{N(\alpha)}{1 - \alpha} \int_{0}^{t} \sin\left[\frac{a}{1-a}(t-\tau)\right]\frac{d}{d\tau} f(\tau) \, d\tau\right\},
\]
or, equivalently, if $*$ is convolution,
\[
\mathcal{L}\left\{\mathfrak{D}_{\sin}^\alpha f(t)\right\}
=
\frac{N(\alpha)}{1 - \alpha}\mathcal{L}\left\{\sin\left(\frac{a}{1-a}t\right)*\frac{d}{dt} f(t)\right\},
\]
or, equivalently, 
\[
\mathcal{L}\left\{\mathfrak{D}_{\sin}^\alpha f(t)\right\}
=
\frac{N(\alpha)}{1 - \alpha}\mathcal{L}\left\{\sin\left(\frac{a}{1-a}t\right)\right\}\mathcal{L}\left\{\frac{d}{dt} f(t)\right\},
\]
or, equivalently, 
\[
\mathcal{L}\left\{\mathfrak{D}_{\sin}^\alpha f(t)\right\}
=
\frac{N(\alpha)}{1 - \alpha}\frac{\frac{a}{1-a}}{s^2 + \left(\frac{a}{1-a}\right)^2}
(sF(s)-f(0)),
\]
and consequently we arrive at \eqref{eq6}. We apply now the Laplace transform into \eqref{eq5} and get:
\[
\mathcal{L}\left\{\mathfrak{D}_{\cos}^\alpha f(t)\right\}
=
\mathcal{L}\left\{ \frac{N(\alpha)}{1 - \alpha} \int_{0}^{t} \cos\left[\frac{a}{1-a}(t-\tau)\right]\frac{d}{d\tau} f(\tau) \, d\tau\right\},
\]
or, equivalently, if $*$ is convolution,
\[
\mathcal{L}\left\{\mathfrak{D}_{\cos}^\alpha f(t)\right\}
=
\frac{N(\alpha)}{1 - \alpha}\mathcal{L}\left\{\cos\left(\frac{a}{1-a}t\right)*\frac{d}{dt} f(t)\right\},
\]
or, equivalently, 
\[
\mathcal{L}\left\{\mathfrak{D}_{\cos}^\alpha f(t)\right\}
=
\frac{N(\alpha)}{1 - \alpha}\mathcal{L}\left\{\cos\left(\frac{a}{1-a}t\right)\right\}\mathcal{L}\left\{\frac{d}{dt} f(t)\right\},
\]
or, equivalently, 
\[
\mathcal{L}\left\{\mathfrak{D}_{\cos}^\alpha f(t)\right\}
=
\frac{N(\alpha)}{1 - \alpha}\frac{s}{s^2 + \left(\frac{a}{1-a}\right)^2}
(sF(s)-f(0)),
\]
and consequently we arrive at \eqref{eq7}. The proof is complete.
\\\\
We now consider the (DS) fractional operator~\eqref{eq4} and the (DC) fractional operator~\eqref{eq5}, with $\alpha = \frac{2}{3}$ and $N(\alpha) = 1 - \alpha$. 
Figures~\ref{Figure 1}-\ref{Figure 4} show the results of applying these operators to various functions. In each plot, the horizontal axis represents time \( t \), and the vertical axis represents the value of the corresponding operator applied to the function \( f(t) \). Specifically,  
Figure~\ref{Figure 1} shows the result for \( f(t) = \frac{1}{3}t^3 \),  
Figure~\ref{Figure 2} for \( f(t) = -\cos(t) \),  
Figure~\ref{Figure 3} for \( f(t) = \sin(t) \),  
and Figure~\ref{Figure 4} for \( f(t) = e^t \).

Following this, we apply the (DS) fractional operator~\eqref{eq4}, the (DC) fractional operator~\eqref{eq5}, and the Caputo fractional derivative (C)~\eqref{eq1} with $\alpha = \frac{2}{3}$ to the same set of functions. The results are plotted together for comparison. In each figure, the horizontal axis represents time \( t \), and the vertical axis shows the value of the corresponding operator applied to the function \( f(t) \). 
Figure~\ref{Figure 5} presents the results for \( f(t) = \frac{1}{3}t^3 \) (left), using \( N\left(\frac{2}{3}\right) = 0.9 \) for (DS) and \( N\left(\frac{2}{3}\right) = 4 \) for (DC), and for \( f(t) = -\cos(t) \) (right), with \( N\left(\frac{2}{3}\right) = 0.4 \) for (DS) and \( N\left(\frac{2}{3}\right) = 0.6 \) for (DC).
Figure~\ref{Figure 6} shows the results for \( f(t) = \sin(t) \) (left), using \( N\left(\frac{2}{3}\right) = 0.3 \) for (DS) and \( N\left(\frac{2}{3}\right) = 0.5 \) for (DC), and for \( f(t) = e^t \) (right), with \( N\left(\frac{2}{3}\right) = 0.8 \) for (DS) and \( N\left(\frac{2}{3}\right) = 1.6 \) for (DC).

From the results, we observe that the (DS) fractional operator produces values that are more similar to those of the Caputo (C) fractional derivative than the (DC) operator does. This means that (DS) behaves more like the classical Caputo derivative, while (DC) shows a different behavior. This observation is important when choosing which operator to use, especially in situations where following the traditional behavior of fractional calculus is essential.


Next, we explore the use of the (DS) fractional operator~\eqref{eq4} in modeling the memristor, a nonlinear element found in electrical circuits. The memristor, short for "memory resistor," is a non-volatile electronic component that controls the flow of current and retains memory of the charge that has previously passed through it, even after power is removed. It is considered the fourth fundamental passive circuit element, alongside the resistor, capacitor, and inductor. The concept was first introduced by Leon Chua in 1971, who identified it as the missing link between charge and magnetic flux. The standard model for the memristor is given by the nonlinear relationship:
\[
V(t) = M(q(t))\cdot I(t),
\]
where \(V\) is the voltage across the memristor, \(I\) is the current, and \(M(q)\) is the memristance, which depends on the total electric charge \(q\) that has flowed through the device. Figure~\ref{Figure 7} (left) shows the characteristic \(V\)-\(I\) graph of the memristor, which forms a pinched hysteresis loop, a signature feature of memristive behavior. By applying the (DS) operator to this model, we aim to reformulate the nonlinear expression into a linear form, offering a new perspective on the dynamics of the memristor. Let \( M(q) = q \) and \( I(t) = \sin(t) \). Then the charge is given by \( q(t) = \int_0^t \sin(\tau) \, d\tau = 1 - \cos(t) \), and thus:
\[
V(t) = M(q(t)) \cdot I(t) = [1 - \cos(t)] \sin(t) = \sin(t) - \frac{1}{2} \sin(2t).
\]
If we apply the Laplace transform \( \mathcal{L}\{\cdot\} \), we obtain:
\[
V_{\mathcal{L}}(s) = \frac{1}{s^2 + 1} - \frac{1}{s^2 + 4}=\frac{3}{(s^2 + 1)(s^2 + 4)}= \frac{3}{2} \cdot \frac{2}{s^2 + 4} \cdot \frac{1}{s^2 + 1},
\]
where \( V_{\mathcal{L}}(s) = \mathcal{L}\{V(t)\} \). Thus, by the convolution theorem, the inverse Laplace transform is:
\[
V(t) = \frac{3}{2} \int_0^t \sin(2(t - \tau)) \sin(\tau) \, d\tau.
\]
Hence, if we use the operator~\eqref{eq4} for \(\alpha = \frac{2}{3}\) and \(N\left(\frac{2}{3}\right) = \frac{1}{2}\), the memristor model can be reformulated in the form:
\begin{equation}\label{eq8}
V(t) = \mathfrak{D}_{\sin}^{\frac{2}{3}} q(t)
\end{equation}
Note that the function \(\frac{3}{10} \sin(2t)\) can be approximated by \(t^{-\frac{8}{9}}\) on the interval \(t \in [3.7, 4.24]\), and by \(t^{-\frac{43}{70}}\) on the interval \(t \in [6.82, 7.35]\); see Figure~\ref{Figure 7} (right). Taking into account that \(\frac{d}{dt}q(t) = \sin(t)\), we observe that for \(t \in [3.7, 4.24]\),
\[
\frac{10}{3} \int_0^t \frac{3}{10} \sin[2(t - \tau)] \sin(\tau) \, d\tau 
\approx \frac{10}{3} \int_0^t (t - \tau)^{-\frac{8}{9}} \frac{d}{d\tau} q(\tau) \, d\tau.
\]
Using the Caputo fractional derivative (C)~\eqref{eq1}, this leads to the approximation:
\[
V(t) \approx c_0 \, \mathcal{D}_c^{\frac{8}{9}} q(t), \quad \text{for } t \in [3.7, 4.24], \quad \text{with } c_0 = \frac{10}{3} \, \Gamma\left(\tfrac{1}{9}\right).
\]
Similarly,
\[
V(t) \approx c_0 \, \mathcal{D}_c^{\frac{43}{70}} q(t), \quad \text{for } t \in [6.82, 7.35], \quad \text{with } c_0 = 5 \, \Gamma\left(\tfrac{27}{70}\right).
\]
These results show that the memristor model can be approximated by a linear form involving the Caputo derivative, but only within specific time intervals. Hence, this representation is not valid globally for all \(t > 0\); its accuracy is limited to narrow windows where the approximation of \(\sin(2t)\) by power-law kernels is sufficiently accurate. In contrast, the reformulation proposed in Equation~\eqref{eq8}, using the (DS) operator, provides an exact linear representation of the memristor model that holds for all \(t \in [0, \infty)\). This demonstrates a key advantage of the (DS) operator over the classical Caputo derivative: it enables a complete and global linearization of the nonlinear memristor dynamics without relying on localized approximations.

{\footnotesize
\noindent
\section*{Conclusions}
In this article, we introduced two new fractional operators with trigonometric kernels, developed to improve the modeling and analysis of nonlinear systems with memory, such as the memristor. 
Future work will focus on further investigating the mathematical properties of these operators, including their stability, invertibility, and potential limitations. 

\section*{Statements and Declarations}

\textbf{Competing Interests:} The author declares that there are no financial or non-financial competing interests related to this work.
\\
\textbf{Funding:} The author received no funding for this research.
\\
\textbf{Data Availability:} No data were generated or used during the course of this study.




}

\end{document}